\def\date{{Sept. 8, 2005}}
\magnification=1100
\baselineskip=14truept
\voffset=.75in  
\hoffset=1truein
\hsize=4.5truein
\vsize=7.75truein
\parindent=.166666in
\pretolerance=500 \tolerance=1000 \brokenpenalty=5000
\def\anote#1#2#3{\smash{\kern#1in{\raise#2in\hbox{#3}}}%
  \nointerlineskip}     
\def\note#1{%
  \hfuzz=50pt%
  \vadjust{%
    \setbox1=\vtop{%
      \hsize 3cm\parindent=0pt\eightrm\baselineskip=9pt%
      \rightskip=4mm plus 4mm\raggedright#1%
      }%
    \hbox{\kern-4cm\smash{\box1}\hfil\par}%
    }%
  \hfuzz=0pt
  }
\def\note#1{\relax}

\newcount\equanumber
\equanumber=0
\newcount\sectionnumber
\sectionnumber=0
\newcount\subsectionnumber
\subsectionnumber=0
\newcount\snumber  
\snumber=0

\def\section#1{%
  \subsectionnumber=0%
  \snumber=0%
  \equanumber=0%
  \advance\sectionnumber by 1%
  \noindent{\bf \the\sectionnumber .~#1~}%
}%
\def\subsection#1{%
  \advance\subsectionnumber by 1%
  \snumber=0%
  \equanumber=0%
  \noindent{\bf \the\sectionnumber .\the\subsectionnumber .~#1~}%
}%
\def\prevs{\the\sectionnumber .\the\snumber }

\long\def\Corollary#1{%
  \global\advance\snumber by 1%
  \bigskip
  \noindent{\bf Corollary~\prevs .}%
  \quad{\it#1}%
}%
\long\def\Lemma#1{%
  \global\advance\snumber by 1%
  \bigskip
  \noindent{\bf Lemma~\prevs .}%
  \quad{\it#1}%
}%
\def\Proof{\noindent{\bf Proof.~}}
\long\def\Proposition#1{%
  \advance\snumber by 1%
  \bigskip
  \noindent{\bf Proposition~\prevs .}%
  \quad{\it#1}%
}%
\long\def\Theorem#1{%
  \advance\snumber by 1%
  \bigskip
  \noindent{\bf Theorem~\prevs .}%
  \quad{\it#1}%
}%
\def\ifundefined#1{\expandafter\ifx\csname#1\endcsname\relax}
\def\labeldef#1{\global\expandafter\edef\csname#1\endcsname{\prevs}}
\def\labelref#1{\expandafter\csname#1\endcsname}
\def\label#1{\ifundefined{#1}\labeldef{#1}\note{$<$#1$>$}\else\labelref{#1}\fi}

\def\preveq{(\the\sectionnumber .\the\equanumber)}
\def\neq{\global\advance\equanumber by 1\eqno{\preveq}}

\def\ifundefined#1{\expandafter\ifx\csname#1\endcsname\relax}

\def\equadef#1{\global\advance\equanumber by 1%
  \global\expandafter\edef\csname#1\endcsname{\preveq}%
  \preveq}

\def\equaref#1{\expandafter\csname#1\endcsname}
 
\def\equa#1{%
  \ifundefined{#1}%
    \equadef{#1}%
  \else\equaref{#1}\fi}

\font\eightrm=cmr8%
\font\eightbf=cmb8%
\font\tenbb=msbm10%
\font\sevenbb=msbm7%
\font\fivebb=msbm5%
\newfam\bbfam%
\textfont\bbfam=\tenbb%
\scriptfont\bbfam=\sevenbb%
\scriptscriptfont\bbfam=\fivebb%
\font\tencmssi=cmssi10%
\font\sevencmssi=cmssi7%
\font\fivecmssi=cmssi5%
\newfam\ssfam%
\textfont\ssfam=\tencmssi%
\scriptfont\ssfam=\sevencmssi%
\scriptscriptfont\ssfam=\fivecmssi%
\def\ssi{\fam\ssfam\tencmssi}%

\font\tenmsam=msam10%
\font\sevenmsam=msam7%
\font\fivemsam=msam5%

\def\bb{\fam\bbfam\tenbb}%

\def\hexdigit#1{\ifnum#1<10 \number#1\else%
  \ifnum#1=10 A\else\ifnum#1=11 B\else\ifnum#1=12 C\else%
  \ifnum#1=13 D\else\ifnum#1=14 E\else\ifnum#1=15 F\fi%
  \fi\fi\fi\fi\fi\fi}
\newfam\msamfam%
\textfont\msamfam=\tenmsam%
\scriptfont\msamfam=\sevenmsam%
\scriptscriptfont\msamfam=\fivemsam%
\mathchardef\leq"3\hexdigit\msamfam 36%
\mathchardef\geq"3\hexdigit\msamfam 3E%

\def\d{\,{\rm d}}
\def\D{{\rm D}}
\long\def\DoNotPrint#1{\relax}

\def\Id{{\rm Id}}
\def\limt{\lim_{t\to\infty}}
\def\L{{\ssi L}}
\def\M{{\ssi M}}
\def\oF{\overline F{}}
\def\oG{\overline G{}}
\def\oH{\overline H{}}
\def\oK{\overline K{}}
\def\qed{~~~{\vrule height .9ex width .8ex depth -.1ex}}
\def\ss{\scriptstyle}
\def\T{\hbox{\tencmssi T}}

\def\II{{\bb I}}
\def\NN{{\bb N}\kern .5pt}
\def\ZZ{{\bb Z}}

\centerline{\bf ASYMPTOTIC EXPANSIONS}
\centerline{\bf FOR INFINITE WEIGHTED CONVOLUTIONS}
\centerline{\bf OF LIGHT SUBEXPONENTIAL DISTRIBUTIONS}

\bigskip
 
\centerline{Ph.\ Barbe and W.P.\ McCormick}
\centerline{CNRS, France, and University of Georgia}
 
{\narrower
\baselineskip=9pt\parindent=0pt\eightrm

\bigskip

{\eightbf Abstract.} We establish some asymptotic
expansions for infinite weighted convolutions of distributions having
light subexponential tails. Examples are presented, some showing that in order
to obtain an expansion with two significant terms, one needs to have a general
way to calculate higher order expansions, due to possible cancellations of
terms. An algebraic methodology is employed to obtain the results.

\bigskip


\bigskip
 
\noindent{\eightbf Keywords:} asymptotic expansion,
convolution, tail area approximation, regular variation, subexponential
distributions, infinite order moving average, weighted sum of random
variables.

}

\bigskip\bigskip

\section{Introduction.}
Subexponential distributions (Chistyakov, 1964), as the name suggests,
have tails which decay at a subexponential rate, and hence, provide
good models for heavy-tailed data. These distributions are defined by
a property on their $n$-fold convolutions which make them useful in
various problems involving sums of independent random variables. To be
more precise, we need to introduce some notation. For a distribution
$F$ we write $\oF=1-F$ its tail function and $F^{\star n}$ its
$n$-fold convolution. Next if $f$ and $g$ are two functions, we write
$f\sim g$ to signify that the ratio $f/g$ tends to $1$ at infinity. A
distribution function $F$ is subexponential if and only if
$\overline{F^{\star 2}}\sim 2\oF$, or, equivalently (Chistyakov, 1964),
if $\overline{F^{\star n}}\sim n\oF$ for all positive integers
$n$. Thus, these distributions have particularly nice analytic
features with respect to convolution, and, as such, form an excellent
class that allows both generality and tractibility in asymptotic
analyses of several stochastic models.

When dealing with problems involving extreme values,
subexponential distributions have two possible limiting extremal behavior;
for instance the Pareto type belong to the domain of max-attraction of 
the Fr\'echet distribution, while the Weibull type belong to the domain
of max-attraction of the Gumbel distribution. Hence, as far as extremal 
behavior is concerned, some subexponential distributions are classified
as light tail, in the sense that their extremes have a limiting behavior
which is of the same type as that for the extremes of a normal distribution 
say.

In this paper, we derive tail area expansions for infinite order weighted sums
of light tailed, yet subexponential, random variables. For Pareto type 
distribution, such expansions were obtained in Barbe and McCormick (2005) 
but their form and the proofs differ sharply from those in the current paper.

Examples of applications of subexponentiality include transient
renewal theory (Teugels, 1975; Embrechts and Goldie, 1982), random
walks (Gr\"ubel, 1985; Veraverbeke, 1977), branching processes
(Chistyakov, 1964; Athreya and Ney, 1972; Chover, Ney, and Waigner,
1972), queueing theory (Pakes, 1975), shot noise (Lebedev, 2002),
infinite divisibility (Embrechts, Goldie, and Veraverbeke, 1979), ruin
theory (Asmussen, 1997; Goldie and Kl\"uppelberg, 1998; Tang and
Tsitsiashvili, 2003), compound sums (Cline, 1987; Embrechts, 1985;
Gr\"ubel, 1987), insurance risk theory (Embrechts, Kl\"uppelberg and
Mikosch, 1997), and heavy-tailed linear processes (Rootz\'en, 1986;
Davis and Resnick, 1988; Geluk and De Vries, 2004; Chen, Ng, Tang,
2005). In these papers, inputs to the processes follow a
subexponential distribution and an asymptotic analysis of an output is
obtained, e.g., claim size distribution and ruin probability, age
distribution and expected number of particles alive, service time
distribution and stationary waiting time distribution, or innovation
distribution and tail area for weighted averages. Mainly, these
analyses are first order asymptotic results.  Some second order
results have been obtained. For compound or subordinated
distributions, that is distributions of the form $\sum_{n\geq 0} p_n
F^{\star n}$ with $F$ subexponential, second order results have been
obtained in Omey and Willekens (1987). For finite order convolutions,
Baltrunas and Omey (1998) prove a second order result.

Throughout this paper, we consider an infinite weighted sum $\sum_{i
\in \ZZ}c_i X_i$, where the $X_i$'s are independent and identically
distributed random variables. We write $F$ the distribution function
of the $X_i$'s and $G$ that of the weighted sum.

Davis and Resnick (1988) consider the tail of $G$ when $F$ is in the
domain of max-attraction of the Gumbel distribution, and belongs to a
class slightly larger than the subexponential one. To state in a
compact form the part of their result which we will use, and because
we will be mostly interested in distributions and not in random
variables, we introduce the following notation. If $X$ has
distribution $F$, we write $\M_cF$ the distribution of $cX$. Thus,
whenever $c$ is positive, $\M_cF=F(\cdot/c)$, while $\M_0F$ is the
distribution function of the point mass at $0$. Since it acts on
functions, we call $\M_c$ a multiplication operator. The distribution
of $\sum_{i \in \ZZ}c_i X_i$ is the infinite weighted convolution
$\star_{i\in\ZZ} \M_{c_i}F$.

In this introduction, we also write $c_{(1)}$ the largest weight, 
and $\nu_+$ its multiplicity in the sequence of weights $(c_i)$.

\bigskip

\noindent{\bf Theorem.} (Davis and Resnick, 1988, Proposition 1.3).\quad{\it
  Let $F$ be a subexponential distribution supported on the
  nonnegative half-line and in the domain of max-attraction of the Gumbel
  distribution. Assume furthermore that the $c_i$'s are nonnegative and that
  $$\sum_{i\in\ZZ} c_i^\delta \ 
    \hbox{is finite for some positive $\delta$ less than $1$.}
    \eqno{\equa{BasicC}}
  $$ 
  Then
  $$
    \oG \sim \nu_+\overline{\M_{c_{(1)}}F} \, .
    \eqno{\equa{DavisResnick}}
  $$
}

We remark that Geluk and De Vries (2004) show a more general version
of the above result; see also Chen, Ng, Tang (2005).

As the above result shows, in this case the largest weight controls
the first order asymptotic equivalence. This stands in contrast to
what happens in the regular variation heavy tail case, where all
weights are on an equal footing and appear in the first order
asymptotic equivalence; see Resnick (1987, p.227). One might propose
that the next term in an expansion would be, say, $k_2
\overline{\M_{c_{(2)}}F}$ where $c_{(2)}$ is the second largest
coefficient and $k_2$ equals the multiplicity of $c_{(2)}$ in the
sequence $(c_i)$. This intuition is correct for a narrow range of 
distributions. In section 2, we
consider three regimes within the light subexponential tail one. In the
heaviest of these regimes, the above intuition provides the right
answer for the expansion; see Theorem 2.3.2. In the other two cases, the
missing element in our discussion makes its appearance, namely,
derivatives of the underlying distribution. In fact, our formalism
will not only provide the second order term, but, more generally,
arbitrary many terms provided the distribution $F$ is smooth enough.

We conclude this introduction by noting sufficient conditions for a
distribution to be subexponential have been obtained by several
authors. Simple sufficient conditions in terms of the hazard function
$-\log \oF$ are given in Pitman (1980), Teugels (1975) and
Kl\"uppelberg (1988). Another sufficient condition on $\oF$ is given
by Goldie (1978). Finally, Goldie and Kl\"uppelberg (1998) and Goldie
and Resnick (1988) give sufficient conditions for a distribution to be
both subexponential and in the domain of max-attraction of the Gumbel
distribution.

\medskip

The paper is organized as follows. The next section contains our main
results. A brief illustration of these results is given in section 3,
mostly to show that the main discussion in Barbe and McCormick (2005)
carries over to the setting of the current paper. The proofs of the main
results are in the last section.

\bigskip

\section{Main results.}
In order to state our asymptotic expansions, we need first to discuss
more precisely the class of distribution functions which we will use,
and some algebraic tools which will allow us to both state and compute
these expansions. This will be done in the next two sub-sections. 
Our main results are presented in the third subsection.

Results for some specific random weights require extra developments and
will be reported elsewhere.

\bigskip

\def\prevs{\the\sectionnumber .\the\subsectionnumber .\the\snumber }
\def\preveq{(\the\sectionnumber .\the\subsectionnumber .\the\equanumber)}

\subsection{Hazard rate and smoothly varying functions.}
Consider for the time being a distribution function $F$ supported on
the nonnegative half-line.
Assuming it exists ultimately, the hazard rate $h=F'/\oF$ will be of
primary importance.  Our analysis will show that how the hazard rate
function compares with $t^{-1} \log t$ asymptotically determines the
form of the asymptotic expansion for tail areas.

Writing the tail function as
$$
  \oF(t)=\oF(t_0)\exp \Bigl(-\int_{t_0}^t h(u)\d u\Bigr) \, ,
$$
we see that if $h(t)\sim\alpha /t$ then $\oF$ is regularly varying
with index $\alpha$. On the other hand, if $\limt h(t)=\alpha$, then
$\oF(t) =e^{-\alpha t (1+o(1))}$ has a behavior close to that of an
exponential distribution. Thus, in order to study tail area behavior
for weighted convolutions in the lighter than regularly varying case
but heavier than exponential case, we are led to consider hazard rates
which satisfy
\setbox1=\vbox{\hsize=3.2in\par
  \noindent$h$ is regularly varying,\hfill\vskip3pt\break\noindent
  $\limt th(t)=+\infty \qquad \hbox{ and }\qquad \limt h(t)=0 \, .$
}
$$
  \lower8pt\box1\eqno{\equa{BasicH}}
$$

All the distribution functions satisfying assumption
{\BasicH} are rapidly varying; indeed, writing $\Id$ for the identity
function on the real line, for any $M$, the hazard rate is ultimately
more than $M/\Id$. Therefore, for $t$ large enough and positive $a$ less
than $1$, we have
$$
  {\oF(t)\over \oF(at)} 
  = \exp\Bigl( -\int_{at}^t h(u) \d u\Bigr) 
  \leq \exp \Bigl( -\int_{at}^t M/u \d u\Bigr)
  = a^M \, ,
$$
showing that $\limt \oF(t)/\oF(at)=0$. It follows that if $0<a<b$ then
$\overline{\M_aF}=o(\overline{\M_bF})$.

The representation of distribution functions in terms of
their hazard rate is closely related to the representation for the
class $\Gamma$ of De Haan (1970); see Bingham, Goldie and Teugels
(1989, \S 3.10). In particular, if {\BasicH} holds, then
$\limt {\oF\bigl( t+x/h(t)\bigr)/ \oF(t)} = e^{-x}$
and $F$ belongs to the domain of max-attraction of the Gumbel distribution
(see Bingham, Goldie, Teugels, 1989, Theorem 8.13.4). Moreover, if the
index of regular variation of the hazard rate is negative, then Pitman's
(1980) criterion implies that $F$ is subexponential. Consequently, if the
hazard rate satisfies {\BasicH} and has a negative index of regular variation,
then the corresponding distribution function satisfies the assumptions of
Proposition 1.3 in Davis and Resnick (1988), and {\DavisResnick} holds
provided {\BasicC} holds as well.

For most linear time series models considered in practice, the sequence
has an infinite order moving average representation (see Brockwell and Davis, 
1987, Theorem 3.1.3). In that setting, $F$ is the innovation distribution,
and it is often unrealistic to assume that it is supported on the nonnegative
half-line. A more common assumption then, is that
of tail balance, asserting that $\oF/\overline{\M_{-1}F}$ has a positive
finite limit at infinity. This assumption is not necessary for our results
to hold, and some possible variations should be clear; but it makes
results easier to state, and so we will use it. We then write $h_-$ the
hazard rate of $\M_{-1}F$ when it exists. We then have the representation for
the lower tail of $F$ given by
$$
  F(-t)=F(-t_0)\exp\Bigl( -\int_{t_0}^t h_-(u)\d u\Bigr) \, ,
$$
for $t$ at least $t_0$. Note that if the distribution is tail balanced,
then $\int_{t_0}^t (h-h_-)(u)\d u$ has a limit at infinity. This implies that
if both $h$ and $h_-$ are regularly varying with index greater than $-1$, 
then they are asymptotically equivalent.

Throughout this paper we adopt, without further mention, the following
convention.

\bigskip

\noindent{\bf Convention.} {\it If a distribution function $F$ is tail
balanced and its hazard rate $h$ satistifies a property, then $h_-$
satisfies the same property. To emphazise this convention, we say that
$F$ is strongly tail balanced.
}

\bigskip

For instance, if the property is '$h$ has a continuous $m$-th derivative', 
we mean that $h_-$ has also a continuous $m$-th derivative.

Assuming that the hazard rate is regularly varying is not always
sufficient to obtain higher order results, in particular because
higher order derivatives of the hazard rate will be involved. This
motivates the following definition, introduced in Barbe and McCormick
(2005). It somewhat extends that of a smoothly varying
function as presented in Bingham, Goldie and Teugels (1989).

\bigskip

\noindent{\bf Definition.} {\it
  A function $h$ is smoothly varying of index $\alpha$ and order $m$ if
  it is ultimately $m$-times continuously differentiable and the $m$-th
  derivative $h^{(m)}$ is regularly varying of index $\alpha-m$. We write
  $SR_{\alpha,m}$ the set of all such functions.
}

\bigskip

Recall that $\Id$ is the identity function on the real line. For any
real number $\alpha$ and any nonnegative integer $m$, we write
$(\alpha)_m$ for $\alpha (\alpha-1) \ldots (\alpha-m+1)$.  Note that
if $h$ is smoothly varying of index $\alpha$ and order $m$, then for
any nonnegative integer $k$ at most $m$, the asymptotic equivalence
$h^{(k)}\sim (\alpha)_k h/\Id^k$ holds.

\bigskip

\subsection{Laplace characters.}
To effectively calculate asymptotic expansions, we need
to recall the algebraic method introduced in Barbe and McCormick
(2005) in their analysis of tail areas for linear and related
processes based on innovations with regularly varying tails.  We
recall the terminology and concepts from that paper for which we
shall have use.  We write $\D$ the derivation
operator, that is, if $g$ is a differentiable function, $\D g$ is its
derivative. We write $\D^0$ the identity operator, and $\D^i=\D
\D^{i-1}$ for any positive integer $i$.  For a distribution function
$F$, we write $\mu_{F,i}$ its $i$-th moment.

Because they will appear in our asymptotic expansions, the basic
algebraic objects that we need are the Laplace characters.

\bigskip

\noindent{\bf Definition.} (Barbe and McCormick, 2005). {\quad \it
   Let $F$ be a distribution function having at least $m$ moments. Its Laplace
  character of order $m$ is the differential operator
  $$ 
    \L_{F,m} = \sum_{0\leq i\leq m} {(-1)^i\over i!}\mu_{F,i} \D^i \, .
  $$
}
 
The algebraic role of these operators derives from their property as
morphisms from the convolution semi-group of distribution functions
with at least $m$ moments to the space of differential operators
modulo the ideal generated by $\D^{m+1}$. More specifically define a
bilinear product on (equivalence classes of) differential operators
modulo the ideal generated by $\D^{m+1}$, denoted $\circ$, by
$$ 
  \D^i \circ \D^j = \cases{\D^{i+j} & if $i+j\leq m$, \cr
                           0        & otherwise. \cr}
$$ 
This is just an explicit expression for composition of differential
operators modulo the ideal generated by $\D^{m+1}$. Then, it is a
simple matter, involving only the definition of convolution and the
binomial formula, to show the useful relation
$$ 
  \L_{K \star H,m} = \L_{K,m} \circ \L_{H,m} \, .
$$
We refer to Barbe and McCormick (2005) for a more extensive discussion
of Laplace characters and their practical use.

\bigskip

\subsection{Asymptotic expansions.}
To write concise formulas, we need to introduce a bit more specialized
notation, $\natural$.  This symbol is used to delete a factor from a
convolution product. More precisely, if $H = K_1\star K_2$ is a
convolution product of distributions $K_1$ and $K_2$, then $H\natural
K_2 = K_1$. In particular $G\natural\M_{c_i}F$ is the distribution of
$\sum_{j\in\ZZ\setminus\{i\}} c_jX_j$.

To allow weights and random variables to be negative and still write
concise statements, we introduce a partial ordering on real numbers,
reflecting some form of tail behavior.  If $a$ and $b$ are two real
numbers, we write $a\prec b$ if $\overline{\M_aF}
=o(\overline{\M_bF})$. If there exists a positive $\epsilon$ such that
$\epsilon \overline{\M_a F}\leq \overline{\M_bF}\leq \epsilon^{-1}
\overline{\M_aF}$ ultimately, then we say that $a$ is equivalent to
$b$ and write $a\equiv b$; provided the denominator does not vanish
ultimately, this means that the ratio $\overline{\M_aF}/\overline{\M_bF}$ 
is ultimately bounded away from $0$ and infinity. We also write $a\preceq b$ 
if and only if $a\prec b$ or $a\equiv b$.  Naturally, we define $a\succ b$ 
as $b\prec a$ and $a \succeq b$ as $b\preceq a$.

Since $\oF$ is rapidly varying, if both $a$ and $b$ are positive, 
then $a\prec b$ if $a<b$, while $a\equiv b$ if and only if $a=b$.
When $a$ and $b$ have different signs, the ordering introduced depends
on how both tails of $F$ compare. In particular, if $F$ vanishes in a 
neighborhood of $-\infty$, then $a\equiv 0$ whenever $a$ is negative,
and $\preceq$ defines a total ordering on the real numbers. Note that
the ordering is also total when $F$ is tail-balanced or vanishes in 
a neighborhood of $-\infty$.

We use the notation $(c_{(i)})_{i\in{\NN}^*}$ to represent the
sequence $(c_i)_{i\in\ZZ}$ in nonincreasing $\prec$ order and without
multiplicity; by the latter, we mean that $c_{(i)}$ is strictly
decreasing for the order, except in the case that it eventually reaches
the minimal element $0$ when it stays identically constant $0$ from
that point onward. For instance, if $c_1=c_2=1$ while
$c_3=1/2$ and $c_4=c_5=1/3$ and all the other $c_i$'s are nonnegative less 
than $1/3$, then $c_{(1)}=1$, while $c_{(2)}=1/2$ and $c_{(3)}=1/3$ and so
on. Note that any sequence $(c_i)_{i\in\ZZ}$ with limit 0
at infinity can be put in nonincreasing order indexed by the positive
integers. We say that an element $c_i$ of the sequence is maximal if it
is equivalent to $c_{(1)}$.

Since we will not restrict the sequence $(c_i)_{i\in\ZZ}$ to be nonnegative,
we need to replace condition \BasicC\ by
$$
  \sum_{i\in\ZZ} |c_i|^\delta\hbox{ is finite for some positive $\delta$
  less than $1$.}
  \eqno{\equa{BasicCC}}
$$

The following convention will be convenient.

\medskip

\noindent{\bf Convention.} {\it We say that the 'standard conditions'
hold if \BasicCC\ holds and either $F$ vanishes on a neighborhood of $-\infty$
or is strongly tail balanced.}

\medskip

Minor changes in our proof show that the results of this paper remain
true if instead of requiring $F$ to vanish in a neighborhood of $-\infty$
we impose that $\int_{-\infty}^0 x^m\d F(x)$ is finite whenever we use the 
Laplace character of order $m$ of $F$ in these results.

\medskip

Our first result covers most light subexponential distributions used in
applications, that is those for which $-\log\oF$ is regularly varying of
index positive and less than $1$. In particular, it includes Weibull
distributions; however, the log-normal distribution is excluded.

\Theorem{%
\label{CaseC}%
  Assume that the standard conditions hold and

  \noindent(i) $h(t)\gg t^{-1}\log t$,

  \noindent(ii) $h$ is smoothly varying of negative index and order $m$,

  \smallskip

  \noindent Then
  $$
    \oG = \sum_{i:c_i\equiv c_{(1)}} \L_{G\natural \M_{c_i} F,m}
    \overline{\M_{c_i} F}
    +o(h^m\overline{\M_{c_{(1)}}F}) \, .
  $$
}

To convey the meaning of this formula, let us consider the case when the
sequence $c_i$ is nonnegative, indexed by the positive integers --- that 
is $c_i$ vanishes if $i$ is negative --- decreasing, with $c_1=1$. Writing 
$H=G\natural F$, that is the distribution of $\sum_{i\geq 2}c_iX_i$, the
formula in Theorem {\CaseC} asserts that for any $k$ at most $m$,
$$
  \oG=\oF -\mu_{H,1}\oF' + {\mu_{H,2}\over 2}\oF'' +\cdots + 
  (-1)^k {\mu_{H,k}\over k!} \oF^{(k)} + o(h^k\oF) \, .
$$
More generally, for a sequence $c_i$ of arbitrary sign with maximal element
$c_{(1)}\equiv 1$, and assuming that $F$ is tail-balanced, Theorem \CaseC\
yields the following. Let $\nu_+$ (respectively $\nu_-$)
be the number of maximal positive (respectively negative) $c_i$'s which
are maximal for the $\prec$ ordering. Since $\oF$ and $\overline{\M_{-1}F}$
are rapidly varying and asymptotically of the same order, the maximal 
positive $c_i$'s, if they exist, are all equal, and equal to minus the 
maximal negative $c_i$'s if they exist. Consequently, Theorem \CaseC\ 
asserts that for $k$ at most $m$,
$$\displaylines{\qquad
  \oG=\nu_+\sum_{0\leq i\leq k}{(-1)^i\over i!}\mu_{G\natural\M_{|c_{(1)}|}F,i}
  \D^i\overline{\M_{|c_{(1)}|}F} 
  \hfill\cr\hfill
  {}+ \nu_-\sum_{0\leq i\leq k}{(-1)^i\over i!}\mu_{G\natural\M_{-|c_{(1)}|}F,i}
  \D^i\overline{\M_{-|c_{(1)}|}F} +o(h^k\oF) \, .\qquad\cr}
$$
Qualitatively, this result is very much in agreement with
Davis and Resnick's result, showing that the largest weight drives the
asymptotic behavior of $\oG$. 

\medskip

For practical purposes, Theorem \CaseC\ is applicable over a broad
range of Weibull-like subexponential distributions. However, its scope
does not include the useful subexponential distribution
given by the log-normal distribution. Moreover, Theorem \CaseC\ stands
in contrast with the situation where $\oF$ is regularly varying, where
each weight contributes to a term of the asymptotic expansion
(Barbe and McCormick, 2005). The purpose of the remainder of this
section is to show the delicate transition between regularly varying
tails and those covered by Theorem {\CaseC}.

Theorem {\CaseC} contains the assumption that the hazard rate is asymptotically
much larger than $\Id^{-1}\log$, excluding distributions with tail
$e^{-(\log t)^a}$ when $a$ is greater than $1$ and at most $2$.
The next result covers this class of distributions with $a$
less than $2$ (note $2$ is not allowed yet!). It shows that
this narrow class of distributions, as far as the asymptotic behavior of
weighted convolution goes, retains some features of regularly varying tails,
yet has rapidly varying tails.

\Theorem{%
\label{CaseA}%
  Assume that the standard conditions and {\BasicH} hold, and

  \noindent(i) $h(t)=o(t^{-1}\log t)$ as $t$ tends to infinity

  \line{\noindent (ii)
  $\limsup_{t\to\infty} th(t)^2/ h\bigl(1/h(t)\bigr) < \infty$,
  \hfill{\rm\equa{ThCaseACondh}}
  }

  \smallskip
  \noindent Then, for any positive integer $m$,
  $$
    \oG
    = \sum_{i:c_i\succeq c_{(m)}} \overline{\M_{c_i}F} 
    +o(\overline{\M_{c_{(m)}}F}) \, .
  $$
}

Setting $\eta(t)=th(t)$, condition {\ThCaseACondh} asserts that the supremum
limit of $\eta(t)/\eta\bigl(1/h(t)\bigr)$ is finite. This is the case if 
$\eta$ is of order $\log^a$.

Note that under the assumptions of Theorem {\CaseA}, the index of
regular variation of the hazard rate is $-1$. In the conclusion of
this theorem, the largest term of the expansion is driven by the
largest coefficient $c_i$, the second order term is driven by the
second largest coefficient, and so on. Although formally the expansion
has a similar form as a first order equivalence in the case of
distributions with regularly varying tails, namely
$\sum_{i\in\ZZ}\overline{\M_{c_i}F}$, here, because the tails are
rapidly varying, nonequivalent weights produce asymptotically distinct
contributions to the expansion.  When $\oF$ is rapidly varying, the
$\overline{\M_{c_i}F}$ are ordered by tail dominance, and this yields
formally the formula presented in Theorem \CaseA.

Our last result yields some very unexpected behavior and fills the gap
between Theorems {\CaseC} and \CaseA.  It covers the intermediate
range where the hazard rate is exactly of order $t^{-1}\log t$ and
therefore $\log \oF$ is exactly of order $-\log^2$. In particular, it
covers a distribution often used in practice: the log-normal.  A close look at
the conclusion of that theorem (particularly examining limiting cases
as $\lambda$ tends to $0$ or infinity) reveals a natural progression
in the terms involved in the expansions in Theorems \CaseC\ and
\CaseA\ with only terms equivalent to $c_{(1)}$ in the former result
and as many as desired in the latter result (viz.\ $\lambda=\infty$
and $\lambda=0$).

In this result $\lfloor\cdot\rfloor$
denotes the integer part, that is the largest integer at most equal to the
argument.

\Theorem{%
\label{CaseB}%
  \hfuzz=2pt
  Assume that the standard conditions hold and

  \noindent (i) $h(t)\sim \lambda t^{-1}\log t$ for some positive $\lambda$,

  \noindent (ii) $h$ is smoothly varying (necessarily of
  index $-1$) and order $m$.

  \smallskip
  
  \noindent Then, for any positive integer $k$ at most $m$,
  $$\displaylines{\quad
    \oG
    =\sum_{i\in\ZZ} \II\{\, c_i\succeq c_{(1)}e^{-k/\lambda}\,\} 
    \L_{G\natural \M_{c_i}F,k+\lfloor \lambda\log (|c_i/c_{(1)}|)\rfloor} 
    \overline{\M_{c_i}F} 
    \hfill\cr\hfill
    {}+ o(h^k\overline{\M_{c_{(1)}}F}) \, .
    \quad\cr}
  $$
}%
\hfuzz=0pt%
With regard to this theorem, we note that the presence of the
integer part serves to indicate that
for small $c_i$, the order of the Laplace character may be
reduced. The indicator function ensures that no negative
order appears as an order to a Laplace character. Since
differentiation and multiplication operators change the asymptotic
order of their argument in different ways, it may happen, for example,
that $\overline{M_{c_{(2)}}F}$ is of higher order than
$\D\overline{\M_{c_{(1)}}F}$ yielding a case where the second largest
coefficient has dominance over the largest coefficient in determining
the second order term. This will be illustrated in the next section.

\bigskip


\section{Examples.}
The purpose of this section is to illustrate briefly our three results, 
and show that, as was done more lengthily in Barbe and McCormick (2005) for
regularly varying tails, the asymptotic scale in which an expansion is written
influences which terms of our results should be kept.

The following four examples illustrate successively Theorems \CaseA,
{\CaseB} and \CaseC. Throughout these examples, we consider a sequence
of weights such that $c_i$ is $0$ if $i$ is nonpositive and the
sequence $c_i$ is nonnegative and strictly decreasing for $i$
positive. We also take $c_1$ to be $1$. Finally, the summability
condition {\BasicCC} is assumed to hold. In the fifth example, we
derive a three terms expansion for the marginal distribution of some
linear processes.

\medskip

\noindent{\it Example 1.} Define the functions
$$
  e_1(t)=\exp(-\log^{3/2}t)
  \qquad\hbox{and}\qquad
  e_2(t)=e_1(t)\exp(-\log^{1/4}t) \, ,
$$
so that $e_2=o(e_1)$ at infinity. Consider a distribution function $F$
such that $\oF=e_1+e_2$ ultimately.  Since $\oF\sim e_1$, Davis and
Resnick (1988) shows that $\oG\sim \oF \sim e_1$, because our sequence
$c_i$ is decreasing with $c_1=1$. What is the second order?

Note that for any positive number $a$ less than $1$,
$$
  \log e_1(t/a)-\log e_2(t) \sim {3\over 2} (\log t)^{1/2}\log a \, .
$$
Consequently, $\M_ae_1=o(e_2)$. Hence, for any $i$ at least $2$, we
have $\overline{\M_{c_i}F}=o(e_2)$. It follows from Theorem {\CaseA}
that the relation $\oF=e_1+e_2$, holding ultimately, shows that $\oF$ 
provides not just a first order equivalence to $\oG$ but actually a 
second order result in that $\oG=e_1+e_2+o(e_2)$.

\medskip

\noindent{\it Example 2.} Let $\theta$ be a positive real number. Take 
$h(t)=2\theta t^{-1}\log t$, so that $\oF=\exp (-\theta\log^2)$ on
$[1,\infty)$. Define $\gamma=\sum_{i\not= 1}c_i$. Note $\gamma$ is nonnegative
and finite by our assumptions. We see that
$\mu_{G\natural \M_{c_1}F,1}=\gamma\mu_{F,1}$. Since $c_1=1$, Theorem {\CaseB}
shows that a two terms expansion is either

\vskip 4pt

\halign{case # :\qquad & $#$\hfill &\quad if\quad $#$\hfill\cr
1 & \oF-\gamma\mu_{F,1}\oF' & \overline{\M_{c_2}F}=o(F')\, ,\cr
\noalign{\vskip 3pt}
2 & \oF+\overline{\M_{c_2}F} & \overline{\M_{c_2}F} \gg F' \, , \cr
\noalign{\vskip 3pt}
3 & \oF+\overline{\M_{c_2}F} - \gamma\mu_{F,1}\oF'
                             & \overline{\M_{c_2}F}\asymp F' \cr
}

\vskip 4pt

\noindent We are in the first case if $c_2\leq e^{-1/2\theta}$ and in case 2 
otherwise. Case 3 cannot occur for this specific distribution. But case 3 may 
occur if, for instance, $\oF(t)=\exp(-\log^2t+\log\log t)$ ultimately!
Note also that for the log-normal, $h(t)=t^{-1}\log t +t^{-1}$,
corresponding to $\theta=1/2$.

\medskip

\noindent{\it Example 3.} Assume again that $c_1=1$ and moreover that $k_1=1$.
We write $\mu_i$ for $\mu_{G\natural\M_{c_1}F,i}$. We suppose that the
assumptions of Theorem {\CaseC} are satisfied for $m=3$. The conclusion of
Theorem {\CaseC} yields a seemingly four terms expansion,
$$
  \oG=\oF -\mu_1\oF'+{\mu_2\over 2} \oF''-{\mu_3\over 6} \oF'''
  + o(h^3\oF) \, .
  \eqno{\equa{ExThree}}
$$
However, this formula hides the true number of significant terms in the
expansion. To clarify this remark, we further investigate the meaning of 
{\ExThree} and calculate
$$\eqalign{
  \oF'&{}=-h\oF \, ,\cr
  \oF''&{}=(-h'+h^2)\oF \, ,\cr
  \oF'''&{}=(-h''+3h'h-h^3)\oF \, .\cr
  }
$$
It is then natural to seek an expansion in the scale containing all the
functions involved in these derivatives, that is
$$
  \oF,\quad h\oF,\quad h^2\oF,\quad h^3\oF,\quad h'\oF,\quad hh'\oF,
  \quad h''\oF
$$
which are not $o(h^3\oF)$. Then, for example, a three terms expansion in that 
scale can be obtained from {\ExThree}and is given by
$$
  \oG=\oF +\mu_1h\oF +{\mu_2\over 2} h^2\oF +o(h^2\oF) \, .
$$
In similar fashion, from {\ExThree}, using the expression for the derivatives
of $\oF$, we calculate that a four terms expansion, where now all terms in
the formulas below are significant, is either
$$
  \oG=\oF +\mu_1h\oF +{\mu_2\over 2} h^2\oF - {\mu_2\over 2} h'\oF + o(h'\oF)
$$
or
$$
  \oG=\oF +\mu_1h\oF +{\mu_2\over 2} h^2\oF - {\mu_3\over 6} h^3\oF 
  + o(h^3\oF) \, ,
$$
or
$$
  \oG=\oF +\mu_1h\oF +{\mu_2\over 2} h^2\oF - {\mu_2\over 2}
  h'\oF - {\mu_3\over 6} h^3\oF + o(h^3\oF) \, ,
$$
the first one occurring if $h^3=o(h')$, the second one occurring if
$h'=o(h^3)$ and the third one occurring if both $h^3$ and $h'$ are of
the same asymptotic order. Note that if the index of regular variation
of $-\log \oF$ is less than $1/2$ then we are in the first case, while
if it is more than $1/2$ we are in the second case. When this index is
equal to $1/2$ we can be in any of the three cases a priori.

\medskip

\noindent{\it Example 4.} The purpose of this example is to show that
some cancellation may occur, analogous to that observed for regularly
varying tails in Barbe and McCormick (2005).

Consider a distribution for which ultimately
$$
  \oF(t)=\exp(-\log^{3/2} t)-\exp\bigl(-\log^{3/2}(2t)\bigr) \, .
$$
Since 
$$
  \log^{3/2}(2t)=\log^{3/2}t+{3\over 2} \log^{1/2}t\log 2 \,
  \bigl(1+o(1)\bigr)\, ,
$$ 
we see that
$$
  \oF(t)\sim \exp(-\log^{3/2} t) \, .
$$
Assume that $c_1=1$ and $c_2=1/2$. Theorem {\CaseA} with $m=2$ yields
$$
  \oG(t) =\oF(t)+\oF(2t)+o\bigl(\oF(2t)\bigr) \, ,
$$
apparently giving a two terms expansion. However, note the
cancellation between the second order term of $\oF(t)$ and the leading
term of $\oF(2t)$, so that the above expansion gives only the one term
expansion
$$
  \oG(t)=\exp (-\log^{3/2}t) +o\bigl(\oF(2t)\bigr) \, .
$$

It is clear from this example that a construction similar to that done
at the end of section 3.2 in Barbe and McCormick (2005) can be made:
regardless whether we apply any of the three theorems of the previous
section, for any fixed $m$, we can find a distribution and a sequence
of weights such that the asymptotic expansion of $\oG$ in the theorem
gives in fact a one term expansion, because of cancellations. Hence,
we see once more that obtaining a second order term for an infinite
convolution is not that of obtaining a second order expansion, but
that of obtaining arbitrarily accurate expansions.

\bigskip

\noindent{\it Example 5.} Consider a stationary linear process with 
innovations having a symmetric distribution $F$. The stationary distribution,
$G$, is that of $\sum_{i\in\ZZ} c_iX_i$ for some proper constants $c_i$. 
Since $F$ is symmetric, the odd moments of $G$ vanish.

We assume that the sequence $(c_i)$ has a unique maximal element in
the $\prec$ ordering, equal to $1$ say.

To calculate the first two even moments, 
for any positive integer $k$, write $C_k^*=\sum_{i:c_i\prec 1} c_i^k$.
Then
$$\eqalign{
  \mu_{G\natural \M_{c_{(1)}}F,2}
  &{}=\mu_{G\natural F,2}=C_2^* \mu_{F,2} \cr
  \mu_{G\natural \M_{c_{(1)}}F,4}
  &{}=\mu_{G\natural F,4}=3({C_2^*}^2-C_4^*)\mu_{F,2}^2 + C_4^* \mu_{F,4} 
   \, . \cr
  }
$$

If $F$ obeys the assumptions of Theorem \CaseC\ with $m=6$ say, we have
$$
  \oG =\oF + {C_2^*\over 2} \mu_{F,2}\oF^{(2)} 
  + {3(C_2^{*2}-C_4^*)\mu_{F,2}^2+C_4^*\mu_{F,4}\over 12} \oF^{(4)}
  + O(\oF^{(6)}) \, .
$$
This formula applies for instance to symmetric Weibull type densities
$\alpha e^{-|x|^\alpha}/2\Gamma(1/\alpha)$. It applies also to densities
which are ultimately equal to $x^\beta e^{-x^\alpha}$; compare with Rootzen
(1986) who obtained first order results.

\bigskip


\section{Proofs.}
The proofs are based on the writing of the convolutions in terms of operators
modelled after those introduced in  Barbe and McCormick (2005). For every 
distribution function $K$, define the operator
$$
  \T_K f(t)=\int_{-\infty}^{t/2} f(t-x)\d K(x) \, .
$$
Then,
$$
  \overline{H\star K}=\T_K\oH+\T_H\oK+\M_2(\oH\,\oK) \, .
  \eqno{\equa{ConvTM}}
$$

Throughout the proofs, we write $F_i$ for $\M_{c_i}F$. Without
loss of generality, we assume that the index set of the sequence
$(c_i)$ is $\NN$ and not $\ZZ$.

Define $G_n=F_1\star\cdots \star F_n$. By induction, formula {\ConvTM} yields
$$\displaylines{\
  \oG_n
  =\T_{G_{n-1}}\oF_n 
  + \sum_{2\leq k\leq n} \T_{F_n}\ldots \T_{F_{n-k+2}}\T_{G_{n-k}} \oF_{n-k+1}
  \hfill\cr\hfill
  {}+\M_2(\oG_{n-1}\,\oF_n) 
  + \sum_{2\leq k\leq n} \T_{F_n}\ldots \T_{F_{n-k+2}} 
  \M_2 (\oG_{n-k}\oF_{n-k+1}) \, .
  \ \equa{GnSum}\cr}
$$

\medskip

\noindent{\bf Note.} Throughout this section, we assume that {\BasicH} holds,
possibly without mentioning it.

\bigskip


\subsection{Some lemmas.}
Preliminary lemmas are presented under not necessarily sharpest
conditions. Use of these lemmas appears mainly in the proof of Theorem
\CaseC.  Our first lemma gives the exact order of derivatives of
$\oF$.

\Lemma{%
\label{DkoF}%
  Assume that $h$ is in $SR_{\alpha-1,m}$ for some $\alpha$ nonnegative
  and less than $1$. If $\alpha$ vanishes, assume further that
  $th(t)$ tends  to infinity at infinity.
  For any nonnegative integer $k$ at most $m$,
  $$
    \oF^{(k)}\sim (-1)^kh^k\oF \, .
  $$%
}%

\Proof Using the Fa\`a di Bruno formula (see Roman, 1980, p.809 for the 
formulation that we are using) $\oF^{(k)}$ equals
$$
  \sum_{1\leq i\leq k} {1\over i!} 
  \sum_{\matrix{\ss n_1+\cdots +n_i=k\cr\noalign{\vskip -3pt}
                \ss n_1,\ldots , n_i\geq 1\cr}}
  {k!\over n_1!\ldots n_i!} h^{(n_1-1)} \ldots h^{(n_i-1)} (-1)^i \oF \, .
$$
Since $h$ is smoothly varying of order $m$ and $k$ is at most $m$, for any
index $n_j$ involved in the sum, $h^{(n_j-1)}\sim (\alpha-1)_{n_j-1}h/
\Id^{n_j-1}$.
As usual in the theory of regular variation, when
$(\alpha-1)_{n_j-1}$ vanishes, such asymptotic equivalence must be read as
$h^{(n_j-1)}=o(h/\Id^{n_j-1})$. Thus,
$$
  h^{(n_1-1)} \ldots h^{(n_i-1)}
  \sim (\alpha-1)_{n_1-1}\ldots (\alpha-1)_{n_i-1} h^i /\Id^{k-i} \, ,
$$
where we use $n_1 + \cdots +n_i = k$. The function $h^i/\Id^{k-i}$,
that is $(\Id h)^i/\Id^k$,
is regularly varying of index $i\alpha-k$.  Among such terms, the one
with the largest order is obtained when $i=k$, even if $\alpha$ vanishes,
for in this case $\Id h$ tends to infinity at infinity. This forces all the
$n_i$'s to be $1$ and yields the result.\hfill$\qed$

\bigskip

The next lemma is a Potter type bound which we have already stated and proved
in section 2.

\Lemma{%
\label{Potter}%
  Let $M$ be an arbitrary positive number. There exists $t_1$ such that for
  any $t$ at least $t_1$ and any $\lambda$ at least $1$,
  $$
  \oF(\lambda t) / \oF(t) \leq \lambda^{-M} \, .
  $$
}

Our next lemma shows that for $t$ large enough, the function 
$\oF (t-\cdot)\oF$ is nonincreasing on some interval $[\, M,t/2\,]$, where
$M$ does not depend on $t$.

\Lemma{%
\label{Chenhua}
  Assume that $h$ is in $SR_{\alpha-1,m}$ for some $m$ at least $1$.
  There exists  some positive $t_1$ and $M$ such that for any $t$ at least
  $t_1$, the function $x\mapsto \oF(t-x)\oF(x)$ is nonincreasing in 
  $[\, M,t/2\,]$.
}

\bigskip

\Proof Let $\delta$ be a positive real number such that $\alpha-1+\delta$ is
negative.
Lemma 2.2.4 in Barbe and McCormick (2005) shows that
there exist $t_1$ and $M$ such that for any $t$ at least $t_1$ and
$x$ in $[\, M,t/2]$,
$$\eqalign{
  {1\over h(x)} {\d \over \d x} \log \bigl( \oF(t-x)\oF(x)\bigr)
  &{}= \Bigl( {h(t-x)\over h(x)}-1\Bigr) \cr
  &{}\leq \Bigl(\Bigl( {t\over x}-1\Bigr)^{\alpha-1+\delta}-1\Bigr) \, .\cr
  }
$$
This upper bound is nonpositive in the given range of $x$, implying the
result.\hfill$\qed$

\bigskip

The following result will be instrumental.

\Lemma{%
\label{StochasticOrderBound}%
  Let $K$ be a cumulative distribution function on the nonnegative half-line,
  such that for some $M$ positive, $\oK\leq M\oF$. For any 
  interval $[\,a,b\,]$ and any nonnegative function $f$ continuous 
  and nondecreasing on $[\,a,b\,]$,
  $$
    \int_{[a,b]} f \d K\leq M \Bigl( \int_{[a,b]} f \d F +f\oF(b) \Bigr) \, .
  $$
}

\noindent We will sometimes use this lemma in its limiting form, as $b$
tends to infinity. In particular, if both $F$ and $K$ are supported 
by the nonnegative
half-line, it implies that $\mu_{K,i}\leq M\mu_{F,i}$ for any positive 
integer $i$.

\bigskip 

\Proof Assume that $f$ is differentiable. An integration by parts yields
$$
  \int_{[a,b]} f \d K
  = f\oK (a)-f\oK(b) + \int_{[a,b]} f'\oK(x) \d x \, .
  \eqno{\equa{EqStoch}}
$$
Hence, the left hand side of {\EqStoch} is at most
$M$ times
$$
  f\oF(a)+\int_{[a,b]}f'\oF(x) \d x \, .
$$
But taking $K$ to be $F$ in {\EqStoch} yields that this sum is at most
$$
  \int_{[a,b]}f \d F + f\oF (b) \, .
$$
This prove the lemma when $f$ is differentiable. When it is only continuous, 
use that the differentiable functions are dense in the continuous ones in 
$(L^1[\,a,b\,],\d K)$.\hfill$\qed$

\bigskip

\noindent Combining the two previous lemmas yields the following asymptotic 
rate of decay on the tail of moment integrals for a class of distributions.

\Lemma{%
\label{TailMoment}%
  Let $K$ be a distribution function such that $\oK=O(\oF)$. 
  Assume that the hazard rate function $h$ of $F$ satisfies {\BasicH}. For any 
  nonnegative integers $j$ and $n$,
  $$
    \int_{t/2}^\infty x^j \d K(x) 
    =o(t^{-n}) \, .
  $$%
}%

\Proof There exists some positive $\epsilon$ such that $\epsilon\oK\leq \oF$
ultimately. Using Lemma \StochasticOrderBound, we have for $t$
large enough
$$\eqalign{
  \epsilon\int_t^\infty x^j \d K(x) 
  &{}\leq \int_t^\infty x^j \d F(x) \cr
  &{}=t^j \oF(t) + j t^j \int_1^\infty s^{j-1}\oF(ts)\d s \cr
  &{}\leq t^j \oF(t) + jt^j \oF(t)\int_1^\infty s^{-2} \d s \, . \cr
  }
$$
This bound is rapidly varying, hence it is $o( t^{-n})$ for 
any $n$.\hfill$\qed$

\bigskip


\subsection{Proof of Theorems {\CaseC} and {\CaseB}.}
In this subsection, we first prove Theorems {\CaseC} and {\CaseB} with the
extra assumption that the weights $c_i$'s are nonnegative. We will 
remove this assumption afterwards.
So, from now on, until stated otherwise, we assume that the weights
are nonnegative.

We start by proving a lemma which allows us to neglect terms arising 
from the multiplication operator in formula \GnSum. 

\Lemma{%
\label{CMNeglect}%
  Assume that the hazard rate is regularly varying of negative index and
  that $\liminf_{t\to\infty} th(t)/\log t$ is positive (possibly infinite). 
  Then, for any $k$ nonnegative, $\overline{\M_2F}^{\, 2} = o(h^k\oF)$.
}

\bigskip

\Proof Using the representation for $\oF$, for any $t$ at least $t_0$,
$$
  {\overline{\M_2F}^{\, 2}\over h^k\oF}(t) 
  =\oF(t_0)\exp
  \Bigl( -\int_{t_0}^{t/2}h(u)\d u + \int_{t/2}^t h(u)\d u - k\log h(t)\Bigr)
  \, .
$$
Let $\epsilon$ be a positive real number. Using that $h$ is regularly 
varying of index $\alpha-1$ for some $\alpha$ less than $1$, we have ultimately
$$\displaylines{\qquad
  -\int_{\epsilon t}^{t/2} h(u)\d u+ \int_{t/2}^th(u)\d u 
  \hfill\cr\hfill
  {}= -th(t) \Bigl( \int_\epsilon^{1/2} v^{\alpha-1}\d v -\int_{1/2}^1 
    v^{\alpha-1} \d v\Bigr) \bigl( 1+o(1)\bigr)
  \qquad\cr}
$$
The right hand side above is less than some $-ath(t)$ for some
positive $a$ by taking $\epsilon$ small enough.  

If $\alpha$ is positive, and less than $1$, then the result is clear
since $\log h(t) = o(t h(t))$.

If $\alpha$ vanishes, $a$ can be taken as large as desired by taking
$\epsilon$ small enough.  By assumption $th(t)$ is ultimately at least
like some positive multiple of $\log t$, implying that $\log h(t)$ grows
at most as a multiple of $\log t$; moreover, $th(t)$ is is ultimately
at least some positive multiple of $\log t$. This proves the lemma.\hfill$\qed$

\bigskip

For any positive real number $M$, we write $B(F,M)$ the set of all distribution
functions $K$ such that both $\oK$ and $\overline{\M_{-1}K}$ are less
than $M\oF$ on the nonnegative half-line.

Our next lemma shows that remainder terms stay negligible under the action
of some $\T_K$ operators.

\Lemma{%
\label{BRemainder}%
  Assume that $h$ is regularly varying with negative index
  of regular variation and $F$ is subexponential.
  Let $K$ be a distribution function in some $B(F,M)$. If $g$ is 
  any function which is $o(h^m\oF)$ then, 
  $$
    \T_Kg=o(h^m\oF) \, .
  $$%
}%

\Proof Let $\epsilon$ be a positive real number. Let $t$ be large enough
so that $|g|\leq \epsilon h^m\oF$ on $[\, t/2,\infty)$. Theorem 1.5.3 in 
Bingham, Goldie and Teugels (1987) implies that $h(t)$ is asymptotically
equivalent to $\sup_{s\geq t} h(s)$.
For $t$ large enough, we obtain for some positive $\rho$
$$\eqalign{
  |\T_Kg(t)|
  &{}\leq \epsilon \int_{-\infty}^{t/2} (h^m\oF)(t-x) \d K(x) \cr
  &{}\leq \epsilon 2^\rho h^m(t) \int_{-\infty}^{t/2}\oF(t-x)\d K(x) \, .\cr
  }
$$
Lemma {\StochasticOrderBound} implies that $\int_0^{t/2}\oF(t-x)\d K(x)$
is less than a constant times $\T_F\oF(t)+\M_2\oF^2(t)$, which is at most 
a constant times $\oF^{*2}(t)$, while $\int_{-\infty}^0\oF(t-x) \d K(x)$
is at most $\oF(t)$. The conclusion follows using subexponentiality 
of $\oF$.\hfill$\qed$

\bigskip

The next lemma is at the heart of both Theorems {\CaseC} and {\CaseB}. We will
use it to approximate of some $\T$ operators by Laplace characters,
when acting on derivatives of $\oF$. In the proof, we will use the absolute 
moments
$$
  |\mu|_{K,i}=\int|x|^i\d K(x)\, .
$$

\Lemma{\label{LApprox}%
  For any fixed integer $p$ at most $m$ and any positive $M$,
  $$
    \lim_{\delta\to 0}\limsup_{t\to\infty}
    \,\,\sup
    \biggl| { \int_{-\delta/h(t)}^{\delta/h(t)} \oF^{(p)}(t-x)\d K(x)
             -\L_{K,m-p} \oF^{(p)}(t)\over h^m\oF(t) } \biggr|
    = 0 \, ,
  $$
  where the supremum is taken over all distribution functions $K$
  belonging to $B(F,M)$.
}

\bigskip

\Proof Taylor's formula with remainder term asserts that
$$\displaylines{\quad
  \oF^{(p)}(t-x)
  =\sum_{0\leq j\leq m-p-1} {(-1)^j\over j!} x^j\oF^{(p+j)}(t)
  \hfill\cr\hfill
  {}+ (-1)^{m-p}\int_0^x {y^{m-p-1}\over (m-p-1)!} \oF^{(m)}(t-x+y)\d y \, .
  \quad\cr}
$$
Integrating this equality between $-\delta/h(t)$ and $\delta/h(t)$ with 
respect to $\d K(x)$ shows that
$$
  \Bigl| \int_{-\delta/h(t)}^{\delta/h(t)} \oF^{(p)}(t-x)\d K(x)
  - \L_{K,m-p}\oF^{(p)} (t)\Bigr|
$$
is at most
$$\displaylines{
  \sum_{0\leq j\leq m-p} {|\oF^{(p+j)}(t)|\over j!} 
  \int_{[-\delta/h(t),\delta/h(t)]^c} |x|^j\d K(x)
  \hfill\cr\hfill
  {}+\int_{-\delta/h(t)}^{\delta/h(t)} \Bigl|\int_0^x {y^{m-p-1}\over (m-p-1)!}
  \bigl(\oF^{(m)}(t-x+y)-\oF^{(m)}(t)\bigr) \d y\Bigr|\d K(x) \, .
  \cr}
$$
Combining Lemma \DkoF\ and the proof of \TailMoment, we see that for 
any positive $\delta$,
$$
  \limt \sup_{K\in B(F,M)} 
  { |\oF^{(p+j)}(t)| \int_{[-\delta/h(t),\delta/h(t)]^c} |x|^j\d K(x) \over
    h^m\oF(t) }
  = 0 \, .
$$

Let $\epsilon$ be a positive number.
For any fixed $\delta$, Lemma {\DkoF} and regular variation of $h$ show that
uniformly in $z$ in $[\,-1,1\,]$ say, $\oF^{(m)}\bigl(t-\delta z/h(t)\bigr)
/\oF^{(m)}(t)$ tends to $e^{\delta z}$. Therefore, we can find $\delta$ such 
that, ultimately, for any $z$ in $[\,-1,1\,]$,
$$
  1-\epsilon \leq {\oF^{(m)}\bigl(t-\delta z/h(t)\bigr) \over \oF^{(m)}(t)}
  \leq 1+\epsilon \, .
$$
Consequently, ultimately and uniformly over $B(F,M)$,
$$\displaylines{\qquad
  \int_{-\delta/h(t)}^{\delta/h(t)} \Bigl|\int_0^x {y^{m-p-1}\over (m-p-1)!}
  \Bigl( {\oF^{(m)}(t-x+y)\over \oF^{(m)}(t)} -1 \Bigr) \d y \Bigr|\d K(x)
  \hfill\cr\hskip 2cm
  {}\leq \epsilon \int_{-\delta/h(t)}^{\delta/h(t)} \Bigl|\int_0^x
    {y^{m-p-1}\over (m-p-1)!} \d y \Bigr|\d K(x)
  \hfill\cr\hskip 2cm
  {}\leq \epsilon |\mu|_{K,m-p} \, . 
  \hfill\cr}
$$
Since $\epsilon$ is arbitrary, Lemma {\LApprox} follows from Lemmas {\DkoF}
and \StochasticOrderBound.\hfill$\qed$

\bigskip

As announced, the preceding lemma yields an approximation of some $\T$ operator
by Laplace characters, which we now state and prove.

\Lemma{\label{ApproxTByL}%
  For any $p$ at most $m$,
  $$
    \limt \sup_{K\in B(F,M)} 
    \Bigl| {(\T_K-\L_{K,m})\oF^{(p)}\over h^m\oF}\Bigr|(t) = 0 \, .
  $$
}

\Proof Given Lemma \LApprox, it suffices to show that for any positive 
$\epsilon$ and any $\delta$ small enough,
$$
  \Bigl|\int_{\delta/h(t)}^{t/2} \oF^{(p)}(t-x) \d K(x) \Bigr|
  \leq \epsilon h^m \oF(t)
  \eqno{\equa{ApproxTByLa}}
$$
and
$$
  \Bigl|\int_{-\infty}^{-\delta/h(t)} \oF^{(p)}(t-x)\d K(x)\Bigr|
  \leq \epsilon h^m \oF(t)
  \eqno{\equa{ApproxTByLaa}}
$$
ultimately and uniformly over $B(F,M)$. Using Lemma {\DkoF} and convergence
of the hazard rate to $0$ at infinity, the left hand side of {\ApproxTByLa}
is ultimately bounded by $\int_{\delta/h(t)}^{t/2} \oF(t-x)\d K(x)$. By Lemma
{\StochasticOrderBound}, this is less than
$$
  M\int_{\delta/h(t)}^{t/2} \oF(t-x)\d F(x) + M\oF^2(t/2) \, .
$$
For $t$ large enough, we substitute $\d F(x)$ by $F'(x)\d x$ in this bound. 
Thus, Applying Lemmas \DkoF, \Chenhua\ and \CMNeglect, this upper bound 
is ultimately at most
$$
  Mt\oF\bigl(t-\delta/ h(t)\bigr)\oF\bigr(\delta/ h(t)\bigl) 
  + o\bigl(h^m\oF(t)\bigr) \, .
$$
This implies \ApproxTByLa\ because on one hand,
$\oF\bigl(t-\delta/h(t)\bigr)\sim \oF(t)e^{-\delta}$, and, on the other hand, 
$\oF$ being rapidly varying, $\oF\bigl(\delta/h\bigr)=o(h^p)$ for any 
positive $p$, and, in particular, $\Id \oF(\delta/h)=o(h^m)$.

\DoNotPrint{
Another use of Lemma {\DkoF} allows us to substitute $\d F(x)$ 
by $h(x)\oF(x)\d x$ in that bound, up to doubling $M$ say. Therefore, 
a sufficient condition for {\ApproxTByLa} to hold true is
$$
  \int_{\delta/h(t)}^{t/2} \oF(t-x) \oF (x) \d x +  \oF^2(t/2)
  = o\bigl(h^m\oF(t)\bigr) \, .
  \eqno{\equa{ApproxTByLb}}
$$
To prove this assertion, observe that
$$
  {\oF(t-x)\oF(x)\over \oF(t)\oF(t_0)}
  = \exp \Bigl( -\int_{t_0}^x h(u) \d u + \int_{t-x}^t h(u) \d u\Bigr) \, .
  \eqno{\equa{ApproxTByLc}}
$$
Let $\eta$ be a positive real number less than $1/2$. Uniformly in $x$ between
$\eta t$ and $t/2$,
$$\displaylines{\qquad
  -\int_{\eta x}^x h(u) \d u + \int_{t-x}^t h(u)\d u
  \hfill\cr\noalign{\vskip 3pt}\hskip 2cm
  {}\sim -th(t) \Bigl( \int_{\eta x/t}^{x/t} v^{\alpha-1}\d v
            -\int_{1-x/t}^1 v^{\alpha-1}\d v \Bigr) 
  \hfill\cr\hskip 2cm
  {}\leq {}-{th(t)\over\alpha} \Bigl(\Bigl({x\over t}\Bigr)^\alpha
   (1-\eta^\alpha)-1+\Bigl(1-{x\over t}\Bigr)^\alpha\Bigr)\, .
  \hfill\equa{ApproxTByLcc}\cr
  }
$$
Set $\zeta=x/t$. On $[\, \eta,1/2\,]$, the function 
$$
  \zeta\mapsto \zeta^\alpha-(\eta/2)^2-1+(1-\zeta)^\alpha
$$ 
reaches its minimum at $\eta$. This minimum, $\eta^\alpha (1-2^{-\alpha})-1
+(1-\eta)^\alpha$, is positive provided that $\eta$ is small enough.
Consequently, taking $\eta$ small enough, the bound \ApproxTByLcc\ is of 
the form $-th(t)$ times some positive number. Therefore, there exists some 
positive $a$ such that, ultimately,
$$
  \int_{\eta t}^{t/2}  {\oF(t-x)\oF(x)\over \oF(t)\oF(t_0)}
  \leq t e^{-th(t) a}
  = \exp \bigl( -th(t) a+\log t\bigr) \, .
  \eqno{\equa{ApproxTByLd}}
$$
Since $th(t)\gg \log t$, the left hand side of {\ApproxTByLd} tends to $0$
at infinity. 

The same argument, taking $x$ to be $t/2$ in {\ApproxTByLc}, shows that
$\oF^2(t/2)=o\bigl(h^m\oF(t)\bigr)$. Therefore, to prove {\ApproxTByLb}, we 
only need to show that for any $\delta$ small enough,
$$
  \int_{\delta/h(t)}^{\eta t} \oF(t-x) \oF(x)\d x
  = o\bigl(h^m\oF (t)\bigr) \, .
  \eqno{\equa{ApproxTByLe}}
$$
On that range we still use {\ApproxTByLc}, and note that since $h$ is regularly
varying of negative index at least $-1$, when $x$ is at most $\eta t$, 
we have the ultimate bounds
$$\eqalign{
  \int_{t-x}^t h(u)\d u
  &{}\leq x h(t) \sup_{0\leq s\leq \eta t} h(t-s)/h(t) \cr
  &{}\leq xh(t) (1+2\eta) \cr
  }
$$
thanks to the uniform convergence theorem for regularly varying
functions (see Bingham, Goldie and Teugels, Theorem 1.2.1) and the
condition that $\alpha$ is nonnegative and at most $1$. Moreover, for
any positive $\theta$ less than $1$,
$$
  -\int_0^x h(u)\d u
  \leq -\int_{\theta x}^x h(u)\d u
  \sim -x h(x){1-\theta^\alpha\over \alpha} \, .
$$
Note this asymptotic equivalence remains valid even when $\alpha = 0$
provided one defines the right hand side above to be the limit as 
$\alpha$ approaches $0$, that is $-\log\theta$. This follows by the 
uniform convergence theorem for regularly varying functions.
Therefore, for $t$ large enough
$$
  {\oF (t-x)\oF(x)\over \oF(t)\oF(t_0)} 
  \leq \exp\Bigl( -x h(x) (1-\eta){1-\theta^\alpha\over\alpha} 
                          + xh(t) (1+2\eta)\Bigr)
  \, .
  \eqno{\equa{TLaplacec}}
$$

Let $\rho$ be a positive number less than $1-\alpha$. For $t$ and $x$ 
sufficiently large, using Potter's bound 
$h(t)\leq (1+2\eta) h(x)(t/x)^{\alpha-1+\rho}$, 
we obtain
$$\displaylines{\quad
  x^m{\oF (t-x)\oF(x)\over \oF(t)\oF(t_0)}
  \hfill\cr\hfill
  \eqalign{
  &{}\leq x^m\exp 
    \biggl( -xh(x) \Bigl( (1-\eta){1-\theta^\alpha\over\alpha}
   - (1+2\eta)^2\Bigl( {x\over t}\Bigr)^{1-\alpha-\rho}\Bigr)\biggr) \cr
  &{}\leq \exp 
    \Bigl( -xh(x) \Bigl( (1-\eta){1-\theta^\alpha\over\alpha} 
   - (1+2\eta)^2\eta^{1-\alpha-\rho}\Bigr)+m\log x\Bigr) \, .\cr
  }%
  \cr}
$$
We can take $\eta$ and $\theta$ small enough so that the upper bound
is of the form $\exp \bigl(-axh(x)+m\log x\bigr)$ for some positive
$a$. Note that when $\alpha$ vanishes, $a$ can be taken arbitrary
large by taking $\theta$ small enough. Since $xh(x)/\log x$ is
ultimately bounded from below by some positive number, this implies
that for any $M$ positive, provided $t$ is large enough
$$
  \int_{\delta/h(t)}^{\eta t} x^m{\oF(t-x)\oF(x)\over \oF(t)} \d x
  \leq \int_{\delta/h(t)}^{\eta t} e^{-(M+1)\log x} \d x
  = O\bigl(h^M(t)\bigr)  \, .
$$
Thus, {\ApproxTByLe} holds, completing the proof of \ApproxTByLa.
}

To prove \ApproxTByLaa, we replace $\oF^{(p)}$ by $h^p\oF$ using Lemma
\DkoF. Since $x$ is negative in the range of integration, we bound 
$h^m(t-x)$ by $\sup_{s\geq t}h^m(s)$, which in turn is equivalent
to $h^m(t)$, and bound $\oF(t-x)$ by $\oF(t)$. Then $\ApproxTByLaa$ 
follows, and this completes the proof of Lemma \ApproxTByL.\hfill$\qed$

\bigskip

Note that in Lemma \ApproxTByL, we could replace $\oF$ by $\overline{\M_aF}$
for any positive $a$, since the hazard rate of $\M_aF$ is $a^{-1}h(a\,\cdot\,)$
and has the same properties as $h$ as far as the proof of the lemma is 
concerned.

\bigskip

We now complete the proofs of both Theorems {\CaseB} and {\CaseC},
assuming positivity of the weights. Without
any loss of generality, we assume that $c_{(1)}$ is $1$. Note that
in Theorem {\CaseB} the hazard rate is regularly varying of index
$-1$. While, in Theorem {\CaseC}, the hazard rate is regularly varying of 
index $\alpha-1$ with $0 \leq \alpha <1$.

Recall formula {\GnSum}. Again, using Davis and Resnick's (1988) Proposition
1.3, $\oG_{n-k}\oF_{n-k+1}=O(\oF^2)$. Thus, Lemma {\CMNeglect}
implies 
$$
  \M_2(\oG_{n-k}\oF_{n-k+1})=o(h^m\oF) \, ,
$$ 
and by induction, using Lemma {\BRemainder},
$$
  \M_2(\oG_{n-1}\oF_n)
  + \sum_{2\leq k\leq n} \T_{F_n}\ldots \T_{F_{n-k+2}}
  \M_2(\oG_{n-k}\oF_{n-k+1})
  = o(h^m\oF) \, .
$$
It follows that under either the assumptions of Theorem {\CaseB} or Theorem
{\CaseC}, 
$$
  \oG_n=\T_{G_{n-1}}\oF_n 
  + \sum_{2\leq k\leq n} \T_{F_n}\ldots \T_{F_{n-k+2}}\T_{G_{n-k}}
  \oF_{n-k+1}+o(h^m\oF) \, .
$$
Again, combining Davis and Resnick's (1988) Proposition 1.3 with
Lemma {\ApproxTByL}, we see that
$$
  \T_{G_{n-k}}\oF_{n-k+1} = \L_{G_{n-k},m}\oF_{n-k+1} + o(h^m\oF) \, .
$$
Any Laplace character applied to some $\oF_i$ yields a linear
combination of derivatives of $\oF_i$. Therefore, using Lemmas {\BRemainder}
and {\ApproxTByL} inductively, we see that
$$\eqalignno{
  \oG_n
  &{}=\sum_{1\leq k\leq n} \L_{F_n,m}\circ \cdots \circ \L_{F_{n-k+2},m}
   \circ \L_{G_{n-k},m}\oF_{n-k+1} + o(h^m\oF) \cr
  &{}=\sum_{1\leq k\leq n} \L_{G_n\natural F_{n-k+1},m}\oF_{n-k+1} +o(h^m\oF)
   \, .
  &\equa{BCGnexp}\cr}
$$
Note in the first sum above the term indexed by $k=1$ is to be read as
$\L_{G_{n-1},m}\oF_{n}$.  Write $H_n$ for the distribution function of
$\star_{i>n}\M_{c_i}F$.  We claim that there exists $n$ large enough
so that $\oH_n=o(h^m\oF)$. Since our sequence $(c_i)$ is
nonincreasing, Davis and Resnick's (1988) Proposition 1.3 implies that
$\oH_n\asymp\oF_{n+1}$ and
$$\eqalignno{
  {\overline{\M_{c_{n+1}}F}\over h^m\oF}(t)
  &{}=\exp\Bigl( -\int_t^{t/c_{n+1}}h(u)\d u - m\log h(t)\Bigr) \cr
  &{}=\exp \Bigl( -th(t) {c_{n+1}^{-\alpha}-1\over\alpha}\bigl(1+o(1)\bigr)
   -m\log h(t)\Bigr) \, . \qquad
  &\equa{BCclaim}\cr
  }
$$
Thus, our claim is now evident. Further, {\BCclaim} shows that under the
assumptions of Theorem \CaseC,
in order that $\oH_n=o(h^m\oF)$, it suffices that $n$ be chosen
large enough so that $c_{n+1} <1.$

Next, note that $G=G_n\star H_n$, so that
$$
  \oG= \T_{H_n}\oG_n + \T_{G_n}\oH_n + \M_2(\oH_n\,\oG_n) \, .
$$
Since $\oH_n=o(h^m\oF)$ and Davis and Resnick's (1988) Proposition 1.3
asserts that $\oG_n\asymp\oF$, Lemmas {\BRemainder} and {\CMNeglect} show that 
$$
  \oG=\T_{H_n}\oG_n+o(h^m\oF) \, .
$$
Thus, using Lemmas {\ApproxTByL}, {\BRemainder} and formula {\BCGnexp}, we 
obtain
$$
  \oG=\sum_{1\leq k\leq n} \L_{G\natural F_k,m}\oF_k + o(h^m\oF) \, .
$$

This proves Theorem {\CaseC} since we can take $n$ to be the first integer
for which $c_{n+1}$ is less than $1$, that is to be $k_1$.

This proves Theorem {\CaseB} as well for the following reason.  First,
recall that under the assumptions of that theorem, $\alpha$ vanishes and in
{\BCclaim} we need to read $(c_{n+1}^{-\alpha} -1)/\alpha$ as
equal to $- \log c_{n+1}$.  Then, since $h(t)\sim\lambda t^{-1}\log t$,
equality {\BCclaim} shows that in order to have $\oF_{n+1}=o(h^m\oF)$, we
must have $(-\lambda\log 1/c_{n+1} +m)$ negative, that is $c_{n+1} <
e^{-m/\lambda}$. Hence we have proved that
$$
  \oG=\sum_{n\geq 1}\II\{\, c_n\geq e^{-m/\lambda}\,\} \L_{G\natural F_n,m}
  \oF_n + o(h^m\oF) \, .
$$
Finally, since $\oF_n^{(i)}$ is of order $h^i\oF_n$, it is of order smaller
than $h^m\oF$ as soon as $c_n<e^{-(m-i)/\lambda}$, that is 
$i>m+\lambda\log c_n$. Thus, the differential operators in the 
Laplace characters only make a nonnegligible contribution for
$0 \leq i \leq m+ \lfloor \lambda \log c_n \rfloor$.
Theorem {\CaseB} follows. 

\bigskip

By an obvious symmetry, Theorems {\CaseC} and {\CaseB} follow
if all the $c_i$'s are negative and $F$ is strongly tail balanced.

\bigskip

To remove the sign restriction on the constants, we write $H$ (respectively
$K$) for the distribution of 
$$
  \sum_{i\in\ZZ\, , \,c_i<0} c_iX_i
  \qquad\bigl(\hbox{respectively}
  \sum_{i\in\ZZ\, , \,c_i>0} c_iX_i\,\,\bigr)\, .
$$
Then $G$ is $H\star K$ and
$$
  \oG= \T_H\oK + \T_K\oH + \M_2(\oH\,\oK) \, .
$$
Now that we have obtained the asymptotic expansions for $\oH$ and $\oK$, we
conclude with nearly identical arguments.

\bigskip


\subsection{Proof of Theorem \CaseA.}
We begin by proving a lemma which allows us to neglect terms arising from the
multiplication operator in formula \ConvTM. Recall that under the assumptions
of Theorem {\CaseA}, the hazard rate is regularly varying of index $-1$.

\Lemma{%
\label{AMNeglect}%
  Assume that $h$ is regularly vaying of index $-1$ and {\BasicH} holds.
  For any positive numbers $a$, $b$, and $c$,
  $$
  \overline{\M_a F}\,\,\overline{\M_b F} =o(\overline{\M_c F}) \, .
  $$%
}%

\Proof Since $x\mapsto \overline{\M_xF}$ is nondecreasing on
$(0,\infty)$, it suffices to prove the lemma when $0<c \leq a \leq
b$. Using the representation of $F$ in terms of its hazard rate, we
have for any positive $\delta$ less than $1$ and for any $t$ at least $t_0$,
$$\eqalign{
  \log \Bigl({\overline{\M_aF}\,\,\overline{\M_bF}\over 
              \overline{\M_cF}}(t)
       \Bigr)
  &{}\leq \int_{t/a}^{t/c} h(u)\d u - \int_{\delta t/b}^{t/b}h(u)\d u
   +\log \oF(t_0) \cr
  &{}= th(t) \bigl( \log (a/c) +\log\delta\bigr) \bigl( 1+o(1)\bigr)\, . \cr
  }
$$
Taking $\delta$ small enough ensures that $\log (\delta a/c)$ is negative.
The result follows since $th(t)$ tends to infinity with $t$.\hfill$\qed$

\medskip

Our next lemma is the essential step in the proof. It shows that some $\T$ 
operators are close to the identity.

\Lemma{%
\label{ATId}%
  Let $M$ be a positive number, and let $K$ be a distribution function 
  in $B(F,M)$. Then, under the assumptions of Theorem {\CaseA}, for any 
  positive real numbers $a$ and $b$,
  $$
    (\T_K-\Id)\overline{\M_aF}= o(\overline{\M_bF})
  $$
  at infinity.
}

\bigskip

\Proof Let $\delta$ be a positive number. Consider first the integral
$$\displaylines{\qquad
  \int_{-\delta/h(t)}^{\delta/h(t)} 
  {\oF\bigl( (t-x)/a\bigr)-\oF(t/a)\over \oF(t/b)} \d K(x)
  \hfill\cr\hfill
  {}= {\oF(t/a)\over\oF(t/b)} \int_{-\delta/h(t)}^{\delta/h(t)} 
  \Bigl({\oF\bigl((t-x)/a\bigr)\over \oF(t/a)} -1\Bigr) \d K(x) \, .
  \quad\equa{ATIdEqa}\cr
  }
$$
If the absolute value of $x$ is at most $\delta/h(t)$, then 
$$\eqalign{
  {\oF\bigl( (t-x)/a\bigr)\over \oF(t/a)} 
  &{}= \exp\Bigl( -\int_{t/a}^{(t-x)/a} h(u)\d u \Bigr) \cr
  &{}= \exp \Bigl( -{t\over a}h(t/a)\int_0^{-x/t} 
   {h\bigl( t(1+v)/a\bigr)\over h(t/a)} \d v \Bigr) \, ,\cr
}
$$
is ultimately between $\exp\bigl(-2h(t/a)|x|/a\bigr)$ 
and $\exp\bigl(2h(t/a)|x|/a\bigr)$

Therefore, if $\delta$ is small enough,
$$\eqalign{
  \Bigl|\int_{-\delta/h(t)}^{\delta/h(t)} 
  \Bigl({\oF\bigl( (t-x)/a\bigr)\over \oF(t/a)} -1\Bigr) \d K(x)\Bigr|
  &{}\leq {4\over a}h(t/a) \int_{-\delta/h(t)}^{\delta/h(t)} |x|\d K(x) \cr
  &{}\leq {4\over a}h(t/a) |\mu|_{K,1} \, . \cr
  }
$$
It follows that the absolute value of the right hand side in {\ATIdEqa} is 
ultimately at most
$$
  {4\oF(t/a)\over a\oF(t/b)} h(t/a) |\mu|_{K,1} \, .
  \eqno{\equa{ATIda}}
$$
Recall that $h$ is regularly varying. So, up to the $4|\mu|_{K,1}/a$ factor, 
this upper bound is
$$\displaylines{\quad
  \exp\Bigl( -\int_{t/b}^{t/a} h(u)\d u + \log h(t)+O(1)\Bigr)
  \hfill\cr\hfill
  {}= \exp\Bigl( -th(t)\log (b/a) \bigl( 1+o(1)\bigr)+\log h(t) +O(1)\Bigr) 
  \, . \quad \cr}
$$
Since $h(t)=o(t^{-1}\log t)$ and $th(t)$ tends to infinity, 
$th(t)=o\bigl(\log h(t)\bigr)$. So {\ATIda} tends to $0$ and so 
does the left hand side in {\ATIdEqa}.

It remains to show that
$$\eqalign{
  \int_{\delta/h(t)}^{t/2} \oF\Bigl( {t-x\over a}\Bigr) \d K(x) 
  &{}= o\bigl(\oF(t/b)\bigr) \, , \cr
  \int_{-\infty}^{-\delta/h(t)} \oF\Bigl( {t-x\over a}\Bigr) \d K(x) 
  &{}= o\bigl(\oF(t/b)\bigr) \, , \cr}
$$
as well as
$$
  \oF(t/a)\oK\bigl(\delta/h(t)\bigr)
  = o\bigl(\oF(t/b)\bigr)\, ,
$$
and
$$
  \oF(t/a)K\bigl(-\delta/h(t)\bigr)
  =o\bigl( \oF(t/b)\bigr) \, .
$$
Because of our assumption that both $\oK$ and $\overline{\M_{-1}K}$ are 
$O(\oF)$, and the function $x\mapsto \oF\bigl((t-x)/a\bigr)$ is 
nondecreasing , these assertions are implied by
$$
  {\oF(t/2a)\oF\bigl(\delta/h(t)\bigr)\over \oF(t/b)} 
  = {\oF(t/2a)\over \oF(t/b)}h(t) \, {\oF\bigl(1/h(t)\bigr)\over h(t)}
  = o(1) \, .
$$
This last estimate holds because {\ATIda} tends to $0$ at infinity,
$\oF$ is rapidly varying and $h$ is regularly varying.\hfill$\qed$

\bigskip

From Proposition 1.3 in Davis and Resnick (1988) and using 
Lemma {\AMNeglect}, we see that for any positive $\delta$,
$$
  \M_2 (\oG_{n-k} \oF_{n-k+1}) 
  \asymp \M_2 (\overline{\M_{c_{1}}F}\,\, \overline{\M_{c_{n-k+1}}F})
  = o(\overline{\M_\delta F}) \, .
$$
But then, using the definition of the $\T$ operators,
$$
  \T_{F_{n-k+2}}\M_2(\oG_{n-k}\oF_{n-k+1})
  \leq \M_2\M_2(\oG_{n-k}\oF_{n-k+1})
  =o(\overline{\M_{2\delta}F}) \, ,
$$
and by induction, taking $\delta$ small enough, we see that
$$
  \M_2(\oG_{n-1}\,\oF_n) + \sum_{2\leq k\leq n} \T_{F_n}\ldots \T_{F_{n-k+2}} 
  \M_2 (\oG_{n-k}\oF_{n-k+1}) = o(\oF_n) \, .
$$
Therefore, we have
$$
  \oG_n
  =\T_{G_{n-1}}\oF_n 
  + \sum_{2\leq k\leq n} \T_{F_n}\ldots \T_{F_{n-k+2}}\T_{G_{n-k}} \oF_{n-k+1}
  + o(\oF_n) \, .
$$
We apply Lemma {\ATId} to obtain 
$$
  \T_{G_{n-k}}\oF_{n-k+1}=\oF_{n-k+1}+o(\oF_n) \, .
$$
Then, another use of Lemma {\ATId} shows that
$$
  \T_{F_{n-k+2}} o(\oF_n)=o(\oF_n) \, ,
$$
and proceeding by induction, we obtain the expansion of $\oG_n$,
$$
  \oG_n=\sum_{1\leq k\leq n} \oF_{n-k+1}+o(\oF_n) \, .
$$
To finish the proof, let $H_n$ denote the distribution function 
$\star_{i>n}F_i$, so that $G=G_n\star H_n$. Then, representation
{\ConvTM} gives
$$
  \oG =\T_{H_n}\oG_n + \T_{G_n}\oH_n + \overline{\M_2G_n} \,\,
  \overline{\M_2H_n} \, .
$$
Up to increasing $n$ we can assume that $c_{n+1}<c_n$. Note that by  
Proposition 1.3 in Davis and Resnick (1988), $\oH_n$ is of order 
$\oF_{n+1}$, which is $o(\oF_n)$. Then, by Lemma {\AMNeglect},
$\overline{\M_2 G_n}\,\,\overline{\M_2 H_n}=o(\oF_n)$. Since
$\oH_n=o(\oF_n)$, Lemma {\ATId} shows that $\T_{G_n}\oH_n=
o(\oF_n)$. Yet another application of Lemma {\ATId} shows that 
$$\eqalign{
  \T_{H_n}\oG_n
  &{}=\sum_{1\leq k\leq n} \T_{H_n}\oF_{n-k+1} +o(\oF_n) \cr
  &{}=\sum_{1\leq k\leq n} \oF_{n-k+1} +o(\oF_n) \, . \cr
  }
$$
This proves Theorem {\CaseA} when the weights are nonnegative.

Removal of the sign restriction is accomplished via the decomposition
of a weighted average into two weighted averages in which each sum has all
weights with the same sign as was done at the end of the proof of Theorems
\CaseC\ and \CaseB.

\bigskip

\noindent{\bf Acknowledgements.} We thank Chenhua Zhang for
shortening our original proof of Lemma \ApproxTByL.

\bigskip


\noindent{\bf References}
\medskip

{\leftskip=\parindent
 \parindent=-\parindent
 \par

S.\ Asmussen (1997). {\sl Ruin Probabilities.} World Scientific.

K.B.\ Athreya, P.E.\ Ney (1972). {\sl Branching Processes.}
Springer.

A.\ Baltrunas, E.\ Omey (1998). The rate of convergence for
subexponential distributions. {\sl Liet.\ Mat.\ Rink.}, 38,
1--18; translation in {\sl Lith.\ Math.\ J.}, 38, 1--14.

Ph.\ Barbe, W.P.\ McCormick (2005). Asymptotic expansions for infinite
weighted convolutions of heavy tail distributions and applications,
{\tt http://www.arxiv.org/abs/math.PR/0412537}, submitted.

N.H.\ Bingham, C.M.\ Goldie, J.L.\ Teugels (1989) {\sl Regular Variation},
2nd ed., Cambridge

P.J.\ Brockwell, R.A.\ Davis (1987). {\sl Time Series: Theory and Methods},
Springer.

Y.\ Chen, K.W.\ Ng, Q.\ Tang (2005). Weighted sums of subexponential random
variables and their maxima, {\sl Adv.\ Appl.\ Probab.}, 37, 510--522.

V.P.\ Chistyakov (1964). A theorem on sums of independent
positive random variables and its applications to branching random
processes, {\sl Theor.\ Probab.\ Appl.}, 9, 640--648.

J.\ Chover, P.\ Ney, S. Waigner (1972). Functions
of probability measures, {\sl J.\ Anal.\ Math.}, 26, 255--302.

D.B.H.\ Cline (1987). Convolutions of distributions with
exponential and subexponential tails. {\sl J.\ Austr.\ Math.\ Soc.\
(A)}, 43, 347--365.

R.A.\ Davis, S.I.\ Resnick (1988). Extremes of moving averages of random 
variables from the domain of attraction of the double exponential distribution,
{\sl Stoch.\ Proc.\ Appl.}, 30, 41--68.

L.\ De Haan (1970). {\sl On Regular Variation and its Application to the
Weak Convergence of the Sample Extremes}, Mathematical Centre Tract, 32, 
Amsterdam.

P.\ Embrechts (1985). Subexponential distribution functions and their 
applications: a review, in {\sl Proceedings of the Seventh 
Conference on Probability Theory, (Bra\c{s}ov, 1982)}, 125--136, VNU Sci.\  
Press, Utrecht, 1985.

P.\ Embrechts, C.M.\ Goldie (1980). On closure and factorization
properties of subexponential and related distributions,  {\sl
J.\ Austr.\ Math.\ Soc.\ (A)}, 29, 243--256.

P.\ Embrechts, C.M.\ Goldie (1982). On convolution tails, {\sl
Stoch.\ Proc.\ Appl.}, 13, 263--278.

P.\ Embrechts, C.M.\ Goldie, N.\ Veraverbeke (1979). Subexponentiality and
infinite divisibility, {\sl Z.\ Wahrsch.\  Verw.\ Geb.}, 49, 335--347.

P.\ Embrechts, C.\ Kluppelberg, T.\ Mikosch (1997). {\sl Modelling
Extremal Events}, Springer.

J.L.\ Geluk, G.C.\ De Vries (2004). Weighted Sums
of Subexponential Random Variables and Asymptotic Dependence Between
Returns on Reinsurance Equities, preprint.

C.M.\ Goldie (1978). Subexponential distributions and
dominated-variation tails, {\sl J.\ Appl.\ Prob.}, 15, 440--442.

C.M.\ Goldie, C.\ Kl\"uppelberg (1998). Subexponential
Distributions. In {\sl A Practical Guide To Heavy Tails, Statistical
Techniques and Applications}, R.\ Adler, R.\ Feldman,
and M.\ Taqqu eds., Birkh\"auser, 435--460.

C.M.\ Goldie, S.I.\ Resnick (1988). Distributions that are both subexponential
and in the domain of attraction of an extreme value distribution, {\sl Adv.\
Appl.\ Probab.}, 20, 706--718.

R.\ Gr\"ubel (1985). Tail behaviour of ladder-height
distributions in random walks, {\sl J.\ Appl.\ Probab.}, 22, 705--709.

R.\ Gr\"ubel (1987). On subordinated distributions and generalized renewal
measures, {\sl Ann.\ Probab.}, 15, 394--415.

C.\ Kl\"uppelberg (1988). Subexponential distributions and
integrated tails, {\sl J.\ Appl.\ Probab.}, 25, 132--141.

A.V.\ Lebedev (2002). Extremes of Subexponential Shot Noise. {\sl
Mathematical Notes}, 71, 206--210.

E.\ Omey (1994). On the difference between the product and the
convolution product of distribution functions. {\sl Publ.\ Inst.\ Math.\
(Beograd) (NS)}, 55, 111--145.

E.\ Omey, E.\ Willekens (1987). Second-order behaviour of
distributions subordinate to a distribution with finite mean. {\sl
Comm.\ Statist.\ Stochastic Models}, 3, 311--342.

A.G.\ Pakes (1975). On the tails of waiting-time distributions,
{\sl J.\ Appl.\ Prob.}, 12, 555--564.

E.J.G.\ Pitman (1980). Subexponential Distribution
Functions, {\sl J.\ Austr.\ Math.\ Soc.\ (A)}, 29, 337--347.

S.\ Roman (1980). The formula of Fa\`a di Bruno, {\sl Amer.\ Math.\ Monthly},
87, 805--809. 

H.\ Rootzen (1986). Extreme value theory for moving average
processes. {\sl Ann.\ Probab.}, 14, 612--652.

Q.\ Tang, G.\ Tsitsiashvili (2003). Randomly weighted sums of
subexponential random variables with application to ruin theory.  {\sl
Extremes}, 6, 171--188.

J.L.\ Teugels (1975). The class of subexponential
distributions, {\sl Ann.\ Probab.}, 3, 1000--1011.

N.\ Veraverbeke (1977). Asymptotic behaviour of Weiner-Hopf
factors of a random walk, {\sl Stoch.\ Proc.\  Appl.}, 5, 27--37.

}

\vskip .5in

\setbox1=\vbox{\halign{#\hfil&\hskip 40pt #\hfill\cr
  Ph.\ Barbe            & W.P.\ McCormick\cr
  90 rue de Vaugirard   & Dept.\ of Statistics \cr
  75006 PARIS           & University of Georgia \cr
  FRANCE                & Athens, GA 30602 \cr
                        & USA \cr
                        & bill@stat.uga.edu \cr}}
\box1

\vfill
\bye